\documentclass[]{amsart}
\usepackage[mathscr]{eucal}
\usepackage{mathrsfs}   
\usepackage{amsmath, amsthm}

\usepackage{setspace}

\newcommand{\Real}{\mathbb{R}}
\newcommand{\Natural}{\mathbb{N}}

\newcommand{\dZ}{\dot{Z}}
\newcommand{\de}{\dot{e}}
\newcommand{\df}{\dot{f}}

\theoremstyle{plain}
\newtheorem{theorem}{Theorem}[section]
\newtheorem{lemma}[theorem]{Lemma}
\newtheorem{prop}[theorem]{Proposition}
\newtheorem{corollary}[theorem]{Corollary}

\theoremstyle{remark}

\newtheorem{remark}[theorem]{Remark}

\numberwithin{equation}{section}

\title{An Algebraic Approach to the Cameron-Martin-Maruyama-Girsanov Formula}
\author{Jir\^o Akahori \and Takafumi Amaba \and Sachiyo Uraguchi}
\thanks{{\bf Mathematics Subject Classification} Primary 60H99; Secondary 60H07, 81T99}
\address{(J. Akahori and T. Amaba) Ritsumeikan University, 
1-1-1 Nojihigashi, Kusatsu, Shiga, 525-8577, Japan}
\address{(S. Uraguchi) Mitsubishi Tokyo UFJ Bank}
\email{(T. Amaba) capca0310@gmail.com} 

\date{}

\begin{document}
\maketitle

\begin{abstract}
In this paper, we will give a new perspective 
to the Cameron-Martin-Maruyama-Girsanov formula
by giving a totally algebraic proof to it.
It is based on the exponentiation of the Malliavin-type
differentiation and its adjointness.   
\end{abstract}
\section{Introduction.}
Let $(\mathscr{W},\mathscr{B}(\mathscr{W}),\gamma )$ 
be the Wiener space on the interval $[0,1]$, that is, 
$\mathscr{W}$ is the set of all continuous 
paths in $\Real$ defined on $[0,1]$ which starts from zero, 
$\mathscr{B}(\mathscr{W})$ is the $\sigma$-field generated 
by the topology of uniform convergence. 
and $\gamma$ is the Wiener measure on the measurable space 
$(\mathscr{W},\mathscr{B}(\mathscr{W}))$. 
Then the canonical Wiener process $(W(t))_{t\geq 0}$ is 
defined by $W(t,w)=w(t)$ for $0\leq t\leq 1$ and $w \in \mathscr{W}$. 

Let $\mathscr{H}$ denote 
the Cameron-Martin subspace of $\mathscr{W}$, 
i.e., $h(t) \in \mathscr{W}$ belongs to $\mathscr{H}$ if and only if 
$h(t)$ is absolutely continuous in $t$ and the derivative 
$\dot{h}(t)$ is  square-integrable. 
Note that $\mathscr{H}$ is a Hilbert space under the inner product
\begin{equation*}
\langle h_{1},h_{2} \rangle _{\mathscr{H}} 
= \int _{0}^{1} \dot{h_{1}}(t)\dot{h_{2}}(t)dt, 
\qquad h_{1},h_{2}\in \mathscr{H}.
\end{equation*}
It is a fundamental fact in stochastic calculus that
the Cameron-Martin (henceforth CM) formula (see, e.g. \cite{Ma}, pp 25) 
in the following form holds:
\begin{equation}\label{Cameron-Martin}
\begin{split}
& \int _{\mathscr{W}} \!\! F(w+\theta) \gamma (dw) \\
& \hspace{1.5cm} 
= \int _{\mathscr{W}} \!\! F(w) \exp \Big\{ \int _{0}^{1} 
\!\! \dot{\theta}(t)dw(t) - \frac{1}{2}\int _{0}^{1} \!\! 
\dot{\theta}(t)^{2}dt \Big\} \gamma (dw) \\
\end{split}
\end{equation}
where $F$ is a bounded measurable function on $\mathscr{W}$ 
and $\theta \in \mathscr{H}$. 

The motivation of the present study comes 
from the following observation(s). 
In the above CM formula (\ref{Cameron-Martin}), 
the integrand of the left-hand-side can be seen 
as an action of a translation operator, 
which is an exponentiation of a differentiation $D_{\theta}$:
\begin{equation}\label{CM2}
\int _{\mathscr{W}} \!\! F(w+\theta) 
\gamma (dw)\ ``\!=" \ E\Big[ e^{D_{\theta}}F \Big].
\end{equation}
On the other hand, the right-hand-side can be seen 
as a ``coupling" of the exponential martingale and $F$:
\begin{align*}
& \int _{\mathscr{W}} \!\! F(w) 
\exp \Big\{ \int _{0}^{1} \!\! \dot{\theta}(t)dw(t) 
- \frac{1}{2}\int _{0}^{1} 
\!\! \dot{\theta}(t)^{2}dt \Big\} \gamma (dw) \\
&\hspace{20mm}= \left \langle  F, 
\exp \Big\{ \int _{0}^{1} \!\! \dot{\theta}(t)dW(t) 
- \frac{1}{2}\int _{0}^{1}
\!\! \dot{\theta}(t)^{2}dt \Big\} \right \rangle .
\end{align*}
Since we can read the right-hand-side of (\ref{CM2}) as
\begin{equation*}
E\Big[ e^{D_{\theta}}F \Big] ``=" \ \left \langle 1,
e^{D_{\theta}}F \right \rangle, 
\end{equation*}
the Cameron Martin formula
\begin{equation*}
\left \langle 1,e^{D_{\theta}}F \right \rangle
``=" \ \left \langle  F, 
\exp \Big\{ \int _{0}^{1} \!\! \dot{\theta}(t)dW(t) 
- \frac{1}{2}\int _{0}^{1} \!\! 
\dot{\theta}(t)^{2}dt \Big\} \right \rangle 
\end{equation*}
leads to the following interpretation:
\begin{equation*}
\exp \Big\{ \int _{0}^{1} \!\! \dot{\theta}(t)dW(t) 
- \frac{1}{2}\int _{0}^{1} \!\! \dot{\theta}(t)^{2}dt \Big\} 
``="\ e^{D_{\theta}^{*}}(1),
\end{equation*}
where $ D^*_\theta $ is an ``adjoint operator" of $ D_\theta $.

The observation, conversely, suggests that 
the CM formula 
could be proved directly by the duality relation between 
$ e^{D_\theta} $ and $ e^{D^*_\theta} $, without  
resorting to the stochastic calculus. 
The program is successfully carried out in section \ref{APCM}.
We may say this program runs by the calculus of functionals
of Wiener integrals.

Along the line, we also give an algebraic
proof of the Maruyama-Girsanov (henceforth MG) formula 
(see e.g. \cite[I\!V.38, Theorem (38.5)]{RW}), 
an extension of the CM formula. 
Note that 
MG formula cannot be written in the quasi-invariant 
form as (\ref{Cameron-Martin}), but in the following way:
\begin{equation}\label{MG-0}
\begin{split}
& \int _{\mathscr{W}} \!\! F(w) \gamma (dw) \\
& = \int _{\mathscr{W}} \!\! F(w-Z(w)) 
\exp \Big\{\int _{0}^{1} 
\!\! \dot{Z}(t,w)dw(t) - \frac{1}{2}\int _{0}^{1} \!\! 
\dot{Z}(t,w)^{2}dt \Big\} \gamma (dw). \\
\end{split}
\end{equation}
Here $ Z : \mathscr{W} \to \mathscr{H} $ is a ``predictable" map
such that
\begin{equation*}
\int _{\mathscr{W}} \exp \Big\{ \int _{0}^{1} 
\!\! \dot{Z}(t,w)dw(t) - \frac{1}{2}\int _{0}^{1} \!\! 
\dot{Z}(t,w)^{2}dt \Big\} \gamma (dw) =1. 
\end{equation*}

In this non-linear situation, infinite dimensional vector fields like
$ X_Z \equiv  Z^{i} D_{e_i}  $\footnote{
Here we use Einstein's convention. },  
where $ \{ e_i \} $ is a basis of $ \mathscr{H} $
and $ Z^i = \langle Z, e_i \rangle_{\mathscr{H}} $, may 
play a role of $ D_\theta $ in the linear case, 
but its exponentiation $ e^{X_Z} $ 
does not make sense anymore.  
Instead, we need to consider ``tensor fields" 
$$ D_{Z}^{\otimes n} = Z^{i_1} \cdots Z^{i_n} D_{e_{i_1}} \cdots D_{e_{i_n}} $$ 
and its formal series 
$$ \sum_{n=0}^\infty \frac{1}{n!} D_{Z}^{\otimes n} 
=: \widetilde{e}^{D_{Z}}. $$
We will show in Proposition \ref{heikou_vector} that 
the operator  $ \widetilde{e}^{D_{Z}} $ is the translation by $ Z $;
$ \widetilde{e}^{D_{Z}} ( f(W) ) = f ( W + Z ) $.
To understand MG formula (\ref{MG-0}) in terms of the translation 
operator $ \widetilde{e}^{D_{Z}} $, 
we additionally introduce another sequence $ \{ L_n \} $ of tensor fields 
(see subsection \ref{sec:girsanov_2} for the definition), 
which has the property (Lemma \ref{cmmgtork}) of 
$$
\sum_{n=1}^\infty \frac{1}{n!} L_n
=
\exp \Big\{\int _{0}^{1} \!\! \dot{Z}(t)dw(t) - \frac{1}{2}\int _{0}^{1} \!\! \dot{Z}^{2}(t)dt \Big\}
(\widetilde{e}^{D_{Z}} -1).
$$
Then, as a corollary to the adjoint formula for 
$ L_n $ (
Theorem \ref{new_operator}), 
MG formula can be obtained (Corollary \ref{c_m_m_g}). 

The proof of key theorem (Theorem \ref{new_operator}), however, 
is not ``algebraic" since it involves the use of It\^o's formula. 
This means, we feel, 
a considerable part of the ``algebraic structure" of 
MG formula is still unrevealed. 
We then try to give a purely algebraic proof (=without resorting 
the results from the stochastic calculus) 
to MG formula 
in section \ref{girsanov} at the cost that 
we only consider the case where 
$\dot{Z}$ is a simple predictable process such as
$$
\dot{Z}=\sum_{i=1}^{N} z_i 1_{(t_i,t_{i+1}]} (t).
$$ 
We will consider a family of vector fields like $ z_i D_i $,
where $ D_i $ is the differentiation in the direction 
of $ \int 1_{(t_i,t_{i+1}]} (t) \,dt $. 
A key ingredient in our (second) algebraic proof of MG formula is
the following semi-commutativity:
\begin{equation}\label{causal}
z_i D_j = D_j z_i  \quad \text{if $ j \geq i $}, 
\end{equation}
which may be understood as ``causality".

Actually, the relation (\ref{causal}) implies that
the Jacobian matrix $ DZ = (D_{e_i} Z_j )_{ij} $, if it is defined, 
is upper triangular. In a coordinate-free language, 
it is nilpotent. Equivalently, $ \mathrm{Tr} (DZ)^n = 0 $
for every $ n $, or $ \mathrm{Tr} \wedge^n DZ = 0 $ for every $ n $. 
Since the statements are coordinate-free(=independent of the 
choice of $ \{ e_i \}$), they can be a characterization of the causality
(=predictability) 
in the infinite dimensional setting as well. 
This observation retrieves the result 
in \cite{ZZ} that Ramer-Kusuoka formula (\cite{Ra},\cite{Ku}) 
is reduced to MG formula when $ DZ $ is nilpotent in this sense. 
The observation also implies that Ramer-Kusuoka formula itself can be 
approached in our algebraic way. 
This program has been successfully carried out in \cite{AA}. 

In the present paper, the domains of the operators 
are basically restricted to ``polynomials" (precise
definition of which will be given soon) in order to 
concentrate on algebraic structures. 
We leave in Appendix 
a lemma and its proof   
to ensure the continuity of the operators and 
hence to have a standard version of CM-MG formula. 

To the best of our knowledge, 
an algebraic proof like ours for CMMG formula have never been 
proposed. 
Though we only treat a simplest one-dimensional Brownian case, 
our method can be applied to more general cases 
if only they have a proper 
action of the infinite dimensional Heisenberg algebra. 
The present study is largely motivated by P. Malliavin's 
way to look at stochastic calculus, which for example 
appears in \cite{Ma} and \cite{MT}, and also 
by some operator calculus often found 
in the quantum fields theory (see e.g. \cite{MJD}).
\section{An Algebraic Proof of the Cameron-Martin Formula.}\label{APCM}
\subsection{Preliminaries.}

For any $h\in \mathscr{H}$, we set
$$
[h](w):=\int _{0}^{1}\dot{h}(t)dw(t),\quad w \in \mathscr{W}.
$$
A function $F:\mathscr{W}\rightarrow \Real$ is 
called a {\it polynomial functional} if there exist an $n\in \Natural$, 
$h_{1},h_{2},\cdots ,h_{n}\in \mathscr{H}$ and a polynomial 
$p(x_{1},x_{2},\cdots ,x_{n})$ of $n$-variables such that
$$
F(w)=p \Big( [h_{1}](w),[h_{2}](w),\cdots ,[h_{n}](w) \Big),
\quad w\in \mathscr{W}.
$$
The set of all polynomial functionals is denoted by $\mathbf{P}$. 
This is an algebra over $\Real$ included densely 
in $L^{p}(\mathscr{W})$ for any $1\leq p< \infty$ 
(see e.g. \cite{IW}, pp 353, Remark 8.2). 

Let $\{ e_{i} \} _{i=1}^{\infty}$ be
an orthonormal basis of $\mathscr{H}$. 
If we set
$$
\xi _{i}(w) := [e_{i}](w) = \int _{0}^{1} 
\dot{e_{i}}(t)dw(t), \qquad i=1,2,\cdots
$$
then $\xi _{1},\xi _{2},\cdots$ are mutually independent 
standard Gaussian random variables. 
Let $H_{n}[\xi]$, $n=1,2,\cdots$ be the $n$-th 
Hermite polynomial in $\xi$ defined by the generating function identity
$$
\exp \Big( \lambda \xi -\frac{\lambda ^{2}}{2} \Big) 
= \sum _{n=0}^{\infty} \frac{\lambda ^{n}}{n!}H_{n}[\xi ],
\quad \lambda \in \Real,
$$
and put
$$
{\mathbf \Lambda}:=\left\{ {\mathbf a}=(a_{i})_{i=1}^{\infty}:
\begin{array}{l}
a_{i} \in \mathbb{Z}^{+}, \\
a_{i}=0 \ \text{except for a finite number of $i$'s}
\end{array} \right\}.
$$
We write ${\mathbf a!}:=\prod _{i=1}^{\infty}a_{i}!$ 
for ${\mathbf a}=(a_{i})_{i=1}^{\infty} \in {\mathbf \Lambda}$. 
We define 
$H_{{\mathbf a}}(w) \in {\mathbf P}$, 
${\mathbf a}\in {\mathbf \Lambda}$ by
$$
H_{{\mathbf a}}(w) 
:= \prod _{i=1}^{\infty}H_{a_{i}}[\xi _{i}(w)],\ \quad 
w\in \mathscr{W}
$$
and then $\{ \frac{1}{\sqrt{{\mathbf a!}}}
H_{{\mathbf a}}:{\mathbf a}\in {\mathbf \Lambda} \}$ 
forms an orthonormal basis of $L^{2}(\mathscr{W})$
(see e.g. \cite{IW}).

For a differentiable function $f$ on $\Real$ 
measured by $N_{1}(d\xi)=\frac{1}{\sqrt{2\pi}}e^{-\xi^{2}/2}d\xi$, 
if we define $\partial$ and $\partial ^{*}$ as
$$
\partial f(\xi)=f^{\prime}(\xi)\ \textrm{and}\ 
\partial ^{*}f(\xi)=-\partial f(\xi)+\xi f(\xi),\quad \xi\in\Real
$$
then $\partial ^{*}$ is adjoint to $\partial$ 
on the differentiable class in $L^{2}(\Real ,N_{1})$. 
We note that the $n$-th Hermite polynomial 
$H_{n}$ can be given by $H_{n}[\xi]= (\partial ^{*n}1)(\xi)$.\\
 
\subsection{Directional differentiations and its exponentials}

For a function $F$ on $\mathscr{W}$ and $\theta \in \mathscr{H}$, 
the differentiation of $F$ in the direction 
$\theta$ $D_{\theta}F$ is defined by
$$
D_{\theta}F(w):=\lim _{\varepsilon \rightarrow 0} 
\frac{1}{\varepsilon}\Big\{ F(w+\varepsilon \theta )-F(w) \Big\},
\quad w\in \mathscr{W}
$$
if it exists(see e.g. \cite{IW}). 
Note that $D_{\theta}F(w)$ is linear 
in $\theta$ and $F$ and satisfies the Leibniz' formula 
$D_{\theta}(FG)(w)=D_{\theta}F(w)\cdot G(w) + F(w)D_{\theta}G(w)$ 
for functions $F$ and $G$ on $\mathscr{W}$ 
such that $D_{\theta}F(w)$ and $D_{\theta}G(w)$ exist. 
If $F(w)$ is of the form $F(w)=f([h](w))$ 
where $f$ is a differentiable function on $\Real$ 
and $h\in \mathscr{H}$, then we have
\begin{eqnarray}
D_{\theta}F(w)
= \langle \theta ,h \rangle _{\mathscr{H}} f^{\prime}([h](w)). \label{eq:diff}
\end{eqnarray}
\quad For $\theta \in \mathscr{H}$, 
we define the {\it exponential} of $D_{\theta}$ by
$$
e^{D_{\theta}}F(w) := 
\sum _{n=0}^{\infty}\frac{1}{n!}D_{\theta}^{n}F(w), 
\quad F\in \mathbf{P}\ \text{and}\ w\in \mathscr{W}
$$
which is actually a finite sum by (\ref{eq:diff}).
\begin{lemma}[]\label{Lem1}
For $F,G\in \mathbf{P}$, we have
\begin{eqnarray}
e^{D_{\theta}}(FG)
=e^{D_{\theta}}(F)\cdot e^{D_{\theta}}(G). \label{eq:hom}
\end{eqnarray}
\end{lemma}
\begin{proof} is a straightforward computation:
\begin{align*}
e^{D_{\theta}}(F)\cdot e^{D_{\theta}}(G)
&= \Big( \sum _{n=0}^{\infty}\frac{1}{n!}D_{\theta}^{n}F \Big) 
\cdot \Big( \sum _{n=0}^{\infty}\frac{1}{n!}D_{\theta}^{n}G \Big) \\
&\hspace{-2cm}= \Big( F+D_{\theta}F 
+ \frac{1}{2!}D_{\theta}^{2}F + \frac{1}{3!}D_{\theta}^{3}F 
+\cdots \Big) \\
& \cdot \Big( G+D_{\theta}G 
+ \frac{1}{2!}D_{\theta}^{2}G 
+ \frac{1}{3!}D_{\theta}^{3}G +\cdots \Big) \\
&\hspace{-2cm}=FG+ \Big\{D_{\theta}F\cdot G 
+ FD_{\theta}G \Big\} \\
&\hspace{-1cm}+\Big\{ \frac{1}{2!}D_{\theta}^{2}F\cdot G 
+D_{\theta}F\cdot D_{\theta}G 
+ F\cdot \frac{1}{2!}D_{\theta}^{2}G \Big\} \\
& \hspace{-2cm} +\Big\{ \frac{1}{3!}D_{\theta}^{3}F\cdot G 
+ \frac{1}{2!}D_{\theta}^{2}F\cdot D_{\theta}G 
+ D_{\theta}F\cdot \frac{1}{2!}D_{\theta}^{2}G 
+ F\cdot \frac{1}{3!}D_{\theta}^{3}G  \Big\} \\
& +\cdots \\
&\hspace{-2cm}=FG + D_{\theta}(FG) 
+ \frac{1}{2!}D_{\theta}^{2}(FG) 
+ \frac{1}{3!}D_{\theta}^{3}(FG) + \cdots =e^{D_{\theta}}(FG).
\end{align*}
\end{proof}

\begin{prop}[]\label{Prop1}
For every $F\in \mathbf{P}$, we have
\begin{eqnarray}
e^{D_{\theta}}F(w) = F( w+ \theta ),\quad w\in \mathscr{W}. \label{eq:trans}
\end{eqnarray}
\end{prop}
\begin{proof}
By Lemma \ref{Lem1}, it suffices to show (\ref{eq:trans}) 
for the functional $F\in \mathbf{P}$ of the form 
$F(w)=f([h](w))$ where $f(x)$ is a polynomial in one-variable 
and $h\in \mathscr{H}$. 
Then using (\ref{eq:diff}), we obtain
\begin{align*}
e^{D_{\theta}}F(w) 
&= \sum _{n=0}^{\infty} \frac{1}{n!}D_{\theta}^{n}f([h](w))
=\sum _{n=0}^{\infty} \frac{1}{n!} 
\langle \theta ,h \rangle _{\mathscr{H}}^{n} f^{(n)}([h](w)) \\
&\hspace{-8mm}= \sum _{n=0}^{\infty} 
\frac{1}{n!} f^{(n)}([h](w))
\Big\{ \Big( [h](w)+ \langle \theta ,h \rangle _{\mathscr{H}} \Big) 
-[h](w) \Big\} ^{n} \\
&\hspace{-8mm}= f\Big( [h](w)
+ \langle \theta ,h \rangle _{\mathscr{H}} \Big)
=F(w+\theta ), 
\end{align*}
where $f^{(n)}(x)$ denotes the $n$-th derivative of $f(x)$.
\end{proof}
\subsection{Formal adjoint operator and its exponential.}
In the analogy of $\partial$ and $\partial ^{*}$ in the previous section, 
we define $D_{\theta}^{*}$, $\theta \in \mathscr{H}$ by
\begin{equation*}
D_{\theta}^{*}F(w):=-D_{\theta}F(w) 
+ \int_{0}^{1}\!\!\dot{\theta} (t)dw(t)\cdot F(w), \quad 
F\in \mathbf{P},w\in \mathscr{W}.
\end{equation*}

Let $\{ e_{i} \} _{i=1}^{\infty}$ be an orthonormal basis of 
$\mathscr{H}$ and put $\xi _{i}(w):=[e_{i}](w)$ 
for $i=1,2,\cdots$. 
Then we have
\begin{lemma}[]\label{Lem2}
It holds that 
$$
E \Big[ D_{\theta}H_{n}[\xi _{k}]\cdot H_{m}[\xi _{l}] \Big] 
= E\Big[ H_{n}[\xi _{k}]D_{\theta}^{*} H_{m}[\xi _{l}] \Big]
$$
for any $k,l,m,n=1,2,\cdots$.
\end{lemma}
\begin{proof}
Since $t\mapsto H_{n}[\int _{0}^{t}e_{k}(s)dw(s)]$ ($n\geq 1$) 
is a martingale with initial value zero, if $k\neq l$ the independence 
of $\xi _{k}$ and $\xi _{l}$ and the formula (\ref{eq:diff}) imply 
that both sides become zero when $n, m\geq 1$. 
If $n=m=0$, it is clear that the left-hand side is zero. 
Then the right-hand side equals to
\begin{align*}
E[ D_{\theta}^{*}1 ] = E[ -D_{\theta}1 
+ \int _{0}^{1}\!\!\dot{\theta}(t)dw(t)] 
=E[ \int _{0}^{1}\!\!\dot{\theta}(t)dw(t) ]=0.
\end{align*}
Hence the case $k=l$ suffices. 
Noting that $\xi _{k}$ is a normal Gaussian random variable, 
we have
\begin{align*}
E \Big[ D_{\theta}H_{n}[\xi _{k}]\cdot H_{m}[\xi _{k}] \Big]
&= \langle \theta ,e_{k} \rangle _{\mathscr{H}} 
E \Big[ H_{n}^{\prime}[\xi _{k}] H_{m}[\xi _{k}] \Big] \\
&\hspace{-30mm}= \langle \theta ,e_{k} \rangle _{\mathscr{H}} 
\int _{-\infty}^{\infty}\!\!\!\! 
\partial H_{n}[\xi ]\cdot H_{m}[\xi ]\gamma _{1}(d\xi ) \\
&\hspace{-30mm}= \langle \theta, e_{k} \rangle _{\mathscr{H}} 
\int _{-\infty}^{\infty}\!\!\!\!H_{n}[\xi ] 
\partial ^{*}H_{m}[\xi ]\gamma _{1}(d\xi ) \\
&\hspace{-30mm}=\langle \theta ,e_{k} 
\rangle _{\mathscr{H}} \int _{-\infty}^{\infty} \!\!\!\! H_{n}[\xi ] 
\Big\{ -H_{m}^{\prime}[\xi ]
+\xi H_{m}[\xi ] \Big\} \gamma _{1}(d\xi) \\
&\hspace{-30mm}=\langle \theta ,e_{k} \rangle _{\mathscr{H}} 
E \Big[ H_{n}[\xi _{k}]\Big\{ -H_{m}^{\prime} [\xi _{k}]
+\xi _{k}H_{m}[\xi _{k}] \Big\} \Big] \\
&\hspace{-30mm}= E \Big[ H_{n}[\xi _{k}]
\Big\{ -D_{\theta}H_{m} [\xi _{k}]+\langle \theta ,e_{k} 
\rangle _{\mathscr{H}}\xi _{k}H_{m}[\xi _{k}] \Big\} \Big] .
\end{align*}
Since $\theta$ can be written as 
$\theta =\sum _{k=1}^{\infty}\langle \theta, e_{k} 
\rangle _{\mathscr{H}}e_{k}$, $\int _{0}^{1} \dot{\theta}(t)dw(t)$ 
admits the $L^{2}$-expansion
$$
\int _{0}^{1}\!\! \dot{\theta}(t)dw(t) 
= \sum _{k=1}^{\infty}\langle \theta ,e_{k} 
\rangle _{\mathscr{H}}\xi _{k}.
$$
Now the independence of $\{ \xi _{i} \} _{i=1}^{\infty}$ shows that
\begin{align*}
E\Big[ H_{n}[\xi _{k}] \int _{0}^{1}\!\!\dot{\theta}(t)dw(t) H_{m}
[\xi _{k}] \Big] = E\Big[ H_{n}[\xi _{k}] \langle \theta ,e_{k} 
\rangle _{\mathscr{H}}\xi _{k} H_{m}[\xi _{k}] \Big].
\end{align*}
\end{proof}

\begin{prop}[]\label{Prop2}
For every $F,G\in \mathbf{P}$, it holds that
\begin{eqnarray}
E[ D_{\theta}F\cdot G ] =E[ FD_{\theta}^{*} G ].
\end{eqnarray}
\end{prop}
\begin{proof}
For fixed $F,G\in \mathbf{P}$, 
there exist a positive integer $n\in \Natural$ 
and an orthonormal system $\{ e_{1}, e_{2}, \cdots , e_{n} \}$ 
in $\mathscr{H}$ and polynomials $f(x_{1},x_{2},\cdots ,x_{n})$ 
and $g(x_{1},x_{2},\cdots ,x_{n})$ of $n$-variables such that
\begin{align*}
\begin{array}{ll}
F(w)=f\Big( [e_{1}](w),[e_{2}](w),\cdots , [e_{n}](w) \Big) & \text{and} \\
G(w)=g\Big( [e_{1}](w),[e_{2}](w),\cdots , [e_{n}](w) \Big). & 
\end{array}
\end{align*}
Extend $\{ e_{1}, e_{2}, \cdots , e_{n} \}$ 
to an orthonormal basis $\{ e_{k} \} _{k=1}^{\infty}$ of $\mathscr{H}$. 
Since the degree of the $n$-th Hermite polynomial is exactly $n$,  
$f$ and $g$ can be written as linear combinations 
of finite products of Hermite polynomials. 
From this fact and by the linearity of 
$D_{\theta}$ and $D_{\theta}^{*}$ and the independence, 
$F$ and $G$ may be assumed without loss of generality 
to be of the form
$$
F(w)= \prod _{i=0}^{p}H_{n_{i}}[\xi _{k_{i}}(w)] 
\quad \text{and} \quad 
G(w)= \prod _{i=0}^{p}H_{m_{i}}[\xi _{k_{i}}(w)].
$$
where $\xi _{k}(w)=[e_{k}](w)$ 
and $k_{1}, k_{2}, \cdots , k_{p}$ are mutually distinct. 
Then, using the Leibniz' rule, 
Lemma \ref{Lem2} and the independence of 
$\xi _{1}, \xi _{2}, \cdots$, we have
\begin{align*}
E[ D_{\theta}F\cdot G ]
&= E\Big[ D_{\theta} \prod_{i=1}^{p}H_{n_{i}}[\xi _{k_{i}}] 
\cdot \prod _{i=1}^{p}H_{m_{i}}[\xi _{k_{i}}] \Big] \\
&\hspace{-20mm}= \sum _{i=1}^{p} E\Big[ D_{\theta} 
H_{n_{i}}[\xi _{k_{i}}] \cdot \prod _{j\neq i}H_{n_{j}}[\xi _{k_{j}}] 
\cdot \prod _{i=1}^{p}H_{m_{i}}[\xi _{k_{i}}] \Big] \\
&\hspace{-20mm}= \sum _{i=1}^{p} 
E\Big[ D_{\theta}H_{n_{i}}[\xi _{k_{i}}]
\cdot H_{m_{i}}[\xi _{k_{i}}] \Big] 
E\Big[ \prod _{j\neq i} 
H_{n_{j}}[\xi _{k_{j}}]H _{m_{j}}[\xi _{k_{j}}] \Big] \\
&\hspace{-20mm}= \sum _{i=1}^{p} 
E\Big[ H_{n_{i}}[\xi _{k_{i}}] 
\Big\{ -D_{\theta}H_{m_{i}}[\xi _{k_{i}}] 
+ \langle e_{k_{i}},\theta 
\rangle _{\mathscr{H}}\xi _{k_{i}} 
H_{m_{i}}[\xi _{k_{i}}] \Big\} \Big] \\
& \times E\Big[ \prod _{j\neq i} 
H_{n_{j}}[\xi _{k_{j}}]H _{m_{j}}[\xi _{k_{j}}] \Big] \\
&\hspace{-20mm}= \sum _{i=1}^{p} 
E\Big[ \prod _{j=1}^{p} H_{n_{j}}[\xi _{k_{j}}] 
\Big\{ -D_{\theta}H_{m_{i}}[\xi _{k_{i}}] 
+ \langle e_{k_{i}},\theta \rangle _{\mathscr{H}}\xi _{k_{i}} 
H_{m_{i}}[\xi _{k_{i}}] \Big\} \prod _{j\neq i} 
H _{m_{j}}[\xi _{k_{j}}] \Big] \\
&\hspace{-20mm}= \sum _{i=1}^{p} 
E \Big[ \prod _{j=1}^{p} H_{n_{j}}[\xi _{k_{j}}] 
\Big( -D_{\theta}H_{m_{i}}[\xi _{k_{i}}] \Big) \Big] \\
&+ E \Big[ \prod _{j=1}^{p} H_{n_{j}}[\xi _{k_{j}}] 
\Big\{ \sum _{i=1}^{p}\langle e_{k_{i}},\theta 
\rangle_{\mathscr{H}}\xi _{k_{i}} \Big\} 
\prod _{j=1}^{p} H _{m_{j}}[\xi _{k_{j}}] \Big]. 
\end{align*}
By the orthogonality of $\xi _{1}, \xi _{2}, \cdots$, 
the last term is equal to
$$
E\Big[ \prod _{j=1}^{p} H_{n_{j}}[\xi _{k_{j}}]\cdot \int _{0}^{1} 
\!\! \dot{\theta}(t)dw(t) \prod _{j=1}^{p} 
H_{m_{j}}[\xi _{k_{j}}] \Big],
$$
which completes the proof. 
\end{proof}

\begin{remark}
Note that $\{ D_{\theta}: \theta \in \mathscr{H} \}$ 
determines a linear operator 
$D:\mathbf{P}\rightarrow \mathbf{P}\otimes \mathscr{H}$ 
such that $\langle DF ,\theta \rangle _{\mathscr{H}}=D_{\theta}F$ 
for each $F\in \mathbf{P}$ and $\theta \in \mathscr{H}$. 
The operator can be extended to an operator 
$D:\mathbf{P}\otimes \mathscr{H} 
\rightarrow \mathbf{P}\otimes \mathscr{H} \otimes \mathscr{H}$ 
by $D(F\otimes \theta ) = DF\otimes \theta$.
This operator is commonly used in Malliavin calculus (see e.g. \cite{IW}).
Its ``adjoint" $D^{*}:\mathbf{P}\otimes \mathscr{H} 
\rightarrow \mathbf{P}$ is defined by
$D^{*}F(w)=-\mathrm{tr}\ DF(w) + [F](w),\ F\in 
\mathbf{P}\otimes \mathscr{H}$. 
Then the ``integration by parts formula";
$$
\int _{\mathscr{W}}\langle DF(w), G(w) 
\rangle _{\mathscr{H}} \gamma (dw)
= \int _{\mathscr{W}}F(w)D^{*}G(w)\gamma (dw)
$$
holds for all $F\in \mathbf{P}$ and 
$G\in \mathbf{P}\otimes \mathscr{H}$ (see e.g. \cite{IW}, pp 361). 
Under these notations, $D_{\theta}^{*}F=D^{*}(F\otimes \theta)$ 
for each $F\in \mathbf{P}$ and 
hence the above adjointness follows immediately
from our result and vice versa.
\end{remark}

Next we define the {\it exponential} 
$e^{D_{\theta}^{*}}$ of $D_{\theta}^{*}$ by the formal series
$$
e^{D_{\theta}^{*}}:=
\sum _{n=0}^{\infty}\frac{1}{n!}D_{\theta}^{*n}.
$$
Let $\{ e_{k} \} _{k=1}^{\infty}$ 
be an orthonormal basis of $\mathscr{H}$ as above.
\begin{theorem}[]\label{Thm1}
For every $\theta \in \mathscr{H}$ such that $|\theta |_{\mathscr{H}}=1$, it holds that
\begin{eqnarray}
D_{\theta}^{*n}1
= H_{n}[\int _{0}^{1} \!\! \dot{\theta}(t)dw(t)] 
\in \mathbf{P},\quad n=0,1,2,\cdots \label{eq:gen}
\end{eqnarray}
and hence $e^{D_{\theta}^{*}}1$ can be defined.
In fact, it is the exponential martingale
\begin{eqnarray}
e^{D_{\theta}^{*}}1(w) = 
\exp \Big\{ \int _{0}^{1}\!\! \dot{\theta}dw(t)-\frac{1}{2} \ \Big\},
\quad w\in \mathscr{W}. \label{eq:exp}
\end{eqnarray}
Furthermore, it holds that
\begin{eqnarray}
E\Big[ e^{D_{\theta}} F \Big] = 
E\Big[ F\cdot e^{D_{\theta}^{*}}1 \Big],
\quad F\in \mathbf{P}. \label{eq:adj}
\end{eqnarray}
\end{theorem}
\begin{proof}
We use the induction on $n$ to prove (\ref{eq:gen}). 
It is clear that
\begin{align*}
D_{\theta}^{*}1(w)=\int _{0}^{1}\!\! \dot{\theta}(t) dw(t) =H_{1}[\int _{0}^{1}\!\! \dot{\theta}(t) dw(t)].
\end{align*}
Suppose that (\ref{eq:gen}) holds for $n$. 
We recall that the Hermite polynomials satisfy the identity
\begin{eqnarray}
H_{n+1}[x]=xH_{n}[x]-nH_{n-1}[x].
\end{eqnarray}
Put $\Theta (w):=\int _{0}^{1}\!\! \dot{\theta}(t) dw(t)$. 
Then, noting that $\langle \theta, \theta \rangle _{\mathscr{H}}=1$ 
and using (\ref{eq:diff}),
\begin{align*}
D_{\theta}^{*(n+1)}1
&= D_{\theta}^{*}H_{n}[\Theta]
 = -D_{\theta}H_{n}[\Theta] + \Theta H_{n}[\Theta ] \\
&\hspace{-3mm}=\Theta H_{n}[\Theta ] -nH_{n-1}[\Theta ] 
= H_{n+1}[\Theta ] .
\end{align*}
Hence (\ref{eq:gen}) holds 
for every $n=0,1,2,\cdots$. 
Then (\ref{eq:exp}) follows immediately 
from (\ref{eq:gen}). 

Finally we shall prove (\ref{eq:adj}). 
By using Proposition \ref{Prop2}, for $F\in \mathbf{P}$ we have
\begin{align*}
E\Big[ e^{D_{\theta}} F \Big]
&= \sum _{n=0}^{\infty} \frac{1}{n!} 
E\Big[ D_{\theta}^{n}F \Big]
 = \sum _{n=0}^{\infty} \frac{1}{n!} 
 E\Big[ F\cdot D_{\theta}^{*n}1 \Big]
 = E\Big[ F\cdot e^{D_{\theta}^{*}}1 \Big] .
\end{align*}
\end{proof}
\begin{corollary}[]\label{Cor1}
For every $\theta \in \mathscr{H}$, it holds that
\begin{equation}
e^{D_{\theta}^{*}}1(w)
= \exp \Big\{ \int _{0}^{1}\!\! \dot{\theta}(t)dw(t)
-\frac{1}{2}\int _{0}^{1} \!\! \dot{\theta} (t)^{2}dt \Big\},
\quad w\in \mathscr{W}. \label{eq:exp2}
\end{equation}
Furthermore, it holds that
\begin{eqnarray}
E\Big[ e^{D_{\theta}} F \Big] 
= E\Big[ F\cdot e^{D_{\theta}^{*}}1 \Big],
\quad F\in \mathbf{P}. \label{eq:adj2}
\end{eqnarray}
\end{corollary}
\begin{proof}
Let $\eta =\theta / |\theta |_{\mathscr{H}}$ and then it follows that
$$
D_{\theta}^{*n}1(w)=|\theta|_{\mathscr{H}}^{n}D_{\eta}^{*n}1(w)
=|\theta|_{\mathscr{H}}^{n}H_{n}[\int _{0}^{1} \!\! 
\dot{\eta}(t)dw(t)]
$$
for $n=0,1,2,\cdots $ and $w\in \mathscr{W}$ 
by Theorem \ref{Thm1}. 
Hence we have
$$
e^{D_{\theta}^{*}}1(w) 
= \sum _{n=0}^{\infty} 
\frac{|\theta |_{\mathscr{H}}^{n}}{n!} 
H_{n}[\int _{0}^{1} \!\! \dot{\eta}(t)dw(t)]
=\exp \Big\{ |\theta |_{\mathscr{H}} \!\! \int _{0}^{1} 
\!\! \dot{\eta}(t)dw(t) 
- \frac{ |\theta |_{\mathscr{H}}^{2} }{2} \Big\} .
$$
The identity (\ref{eq:adj2}) can be shown by the same argument 
as Theorem \ref{Thm1}.
\end{proof}

Now, we have 
the Cameron-Martin formula in this polynomial framework.
\begin{corollary}\label{Cor2}
For every $\theta \in \mathscr{H}$ and $F \in \mathbf{P}$, it holds that
\begin{eqnarray}
\left.
\begin{split}
& &
  \int _{\mathscr{W}}\!\!F(w+\theta )\gamma (dw)
= E \Big[ e^{D_{\theta}}F \Big]
= E \Big[ F\cdot e^{D_{\theta}^{*}}1 \Big] \\
& &
= \int _{\mathscr{W}} \!\! F(w) \exp \Big\{ \int _{0}^{1}\!\! \dot{\theta}dw(t)
  -\frac{1}{2}\int _{0}^{1} \!\! \dot{\theta} (t)^{2}dt \Big\} \gamma (dw). 
\end{split}
\right. \label{eq:CM1}
\end{eqnarray}
\end{corollary}

\section{An Algebraic Proof of MG Formula.}\label{PMG}

In this section, we will give an algebraic proof of 
the MG formula using an adjoint relation similar to (\ref{eq:adj}). 
As we have announced in the introduction, 
for the proof of the adjoint relation 
we will rely on the standard stochastic calculus. 

Let $Z: \mathscr{W} \to \mathscr{H} $ be a predictable map; i.e.
$ \dZ (t)$, $0\leq t\leq 1$ is a predictable process such that
\begin{align*}
\Vert Z \Vert_\mathscr{H}^{2} = \int_{0}^{1} \!\! \dZ(s)^{2}ds < \infty \quad \text{a.s.}
\end{align*}
Suppose $\mathcal{E} (\int \dZ dW )$ is a true martingale 
where for a martingale $M=(M(t))_{0\leq t\leq 1}$ 
the process $\mathcal{E} (M)$ is defined by
\begin{align*}
\mathcal{E} (M)_{t} = \exp \Big\{ M(t)-\frac{1}{2}\langle M \rangle (t)\Big\} .
\end{align*}
\subsection{Infinite dimensional tensor fields.}\label{sec:girsanov_1}

We fix a c.o.n.s. $ \{ e_i: i \in \mathbf{N} \} $ of $ \mathscr{H} $ and 
will write simply $ D_i $ for $ D_{e_i} $ for each $ i \in \mathbf{N} $. 
We define a differentiation along $Z$.
For $\phi \in \mathbf{P}$, we define $D_{Z}$ in the following way:
\begin{equation*}
D_{Z} \phi (W):= \sum _{i=1}^{\infty}\langle Z, e_{i} \rangle (W) D_{i} \phi(W),
\end{equation*} 
where $\langle \cdot, \cdot \rangle$ is the inner product of $\mathscr{H} $. 
Moreover, we define the $n$-th $D_{Z}$, 
which we write as $D_{Z}^{\otimes n}$ by the following:
\begin{equation*}
\begin{split}
D_{Z}^{\otimes n} 
& :=\underbrace{D_{Z} \otimes D_{Z} \otimes \cdots \otimes D_{Z}}_{n} \\
& := \sum _{i, j, k, \cdots} \underbrace{ \langle Z, e_{i} \rangle \langle Z, e_{j} \rangle \langle Z, e_{k} \rangle \cdots }_{n} 
\underbrace{ D_{i} D_{j} D_{k} \cdots }_{n}.
\end{split}
\end{equation*}
Next we define the exponential of $ D_{Z} $ by
the formal series of
\begin{equation*}
\begin{split}
\widetilde{e}^{D_{Z}} 
& := 1 + D_{Z} + \frac {1}{2!}D_{Z}^{\otimes 2}+
\frac{1}{3!}D_{Z} ^{\otimes 3}  + \cdots \\
& = 1 + \sum _{i} \langle Z, e_{i} 
\rangle D_{i} + \frac {1}{2!}\sum _{i, j} 
\langle Z, e_{i} \rangle \langle Z, e_{j} \rangle D_{i} D_{j} \\
& \hspace{2cm} + \frac{1}{3!}\sum _{i, j, k} 
\langle Z, e_{i} \rangle \langle Z, e_{j} \rangle \langle Z, e_{k} 
\rangle D_{i} D_{j} D_{k} + \cdots. \\
\end{split}
\end{equation*}
We denote $\langle Z, e_{i} \rangle $ by $Z_{i}$, 
so we may write $\langle Z, e_{i} \rangle \langle Z, e_{j} \rangle 
D_{i} D_{j}$ as $Z_{i}Z_{j} D_{i} D_{j}$ and furthermore 
$D_{Z}^{\otimes 2}=\sum_{i, j}Z_{i}Z_{j} D_{i} D_{j}$ 
as $\langle Z \otimes Z, \nabla \otimes \nabla \rangle, 
\cdots, 
D_{Z}^{\otimes n}=\langle Z^{\otimes n}, \nabla^{\otimes n} \rangle $, 
and so on. 

\begin{lemma}\label{heikou_e}
For any $k \in \mathbb{N}$, we have 
\begin{equation}\label{heikou_e_f}
\begin{split}
&
\widetilde{e}^{D_{Z}}
\Big(
H_{n_{1}}( \int_{0}^{1} \!\! \de_{m_{1}}dW )
\cdots H_{n_{k}}( \int_{0}^{1} \!\! \de_{m_{k}}dW )
\Big) \\ 
&=
\widetilde{e}^{D_{Z}}
\Big(
H_{n_{1}}( \int_{0}^{1} \!\! \de_{m_{1}}dW )
\Big) 
\cdots \widetilde{e}^{D_{Z}} \Big(
H_{n_{k}}( \int_{0}^{1} \!\! \de_{m_{k}}dW )
\Big).
\end{split}
\end{equation}
\end{lemma}
\begin{proof}
First note that the equation (\ref{heikou_e_f}) is equivalent to 
\begin{equation}\label{heikou_e_sum}
\begin{split}
& \sum_{l=0}^{n_1+\cdots+n_k} \frac{1}{l!}
\langle Z^{\otimes l}, \nabla^{\otimes l} \rangle \Big(
H_{n_{1}}( \int_{0}^{1} \!\! \de_{m_1}dW )
\cdots H_{n_{k}}( \int_{0}^{1} \!\! \de_{m_k}dW )
\Big) \\
&=
\sum_{l_1=0}^{n_1} \frac{1}{l_1!}
\langle Z^{\otimes l_1}, \nabla^{\otimes l_1} \rangle
H_{n_1}( \int_{0}^{1} \!\! \de_{m_1}dW)
\cdots \sum_{l_k=0}^{n_k} \frac{1}{l_k!}
\langle Z^{\otimes l_k}, \nabla^{\otimes l_k} \rangle
H_{n_k}( \int_{0}^{1} \!\! \de_{m_k}dW).\\
\end{split}
\end{equation}
Fixing $l_1, \cdots l_k$ such that $ l_1 \leq n_1, \cdots,l_k \leq n_k $, 
it suffices to prove that the coefficients of 
\[
\nabla^{\otimes l_1}H_{n_1}\nabla^{\otimes l_2}H_{n_2}\cdots \nabla^{\otimes l_k}H_{n_k}
\]
of the left-hand after applying Leibniz rule 
correspond to those of right-hand. 
The coefficients of the left-hand are the following.
\begin{eqnarray*}
\frac {1}{(l_1+l_2+\cdots+l_k)!}
\begin{pmatrix}
l_1+l_2+\cdots+l_k \\
l_1 \\
\end{pmatrix}
\begin{pmatrix}
l_2+\cdots+l_k \\
l_2 \\
\end{pmatrix}
\cdots
\begin{pmatrix}
l_k \\
l_k \\
\end{pmatrix}.
\end{eqnarray*}
This is equal to
$\frac{1}{l_{1}!l_{2}! \cdots l_{k}!}$, so we get (\ref{heikou_e_sum}).
\end{proof}

\begin{prop}\label{heikou_vector}
For $ f \in \mathbf{P}$, we have
\begin{equation}\label{heikou_vector:1}
\widetilde{e}^{D_{Z}} ( f(W) ) = f\left( W + Z \right).
\end{equation}
\end{prop}
\begin{proof}
Since $\widetilde{e}^{D_{Z}}$ is linear and by Lemma \ref{heikou_e}, 
we only prove in the case of $f(W) = H_{n}(\int _{0}^{1} \de_{i}(s)dW_{s})$,
that is, it suffices to show
\begin{align*}
\widetilde{e}^{D_{Z}} \left(
H_{n}( \int _{0}^{1} \!\! \de_{i}(s)dW_{s} )
\right) =
H_{n} \left( \int _{0}^{1} \!\! \de_{i}(s)dW_{s} 
+ \langle Z, e_{i} \rangle
\right).
\end{align*}
By the definition, we have
\begin{align*}
\widetilde{e}^{D_{Z}} \left(
H_{n}( \int _{0}^{1} \!\! \de_{i}(s)dW_{s} )
\right) =
\sum_{k=0}^{n} \binom{n}{k} \langle Z, e_{i} \rangle ^{k}
H_{n-k}( \int _{0}^{1} \!\! \de_{i}(s)dW_{s} ).
\end{align*}
For this, apply $H_{n}(x+y) = \sum_{k=0}^{n}\binom{n}{k}H_{n-k}(x) y^{k}$, then we get (\ref{heikou_vector:1}).
\end{proof}
\subsection{The operator $ L_{n}^Z $.}\label{sec:girsanov_2}

To prove Maruyama-Girsanov formula, we additionally 
introduce a sequence $\{ L_{n}^Z\} $ of new operators associated with $ Z $
as follows. For any $ n \in \mathbb{N}$, $L_{n}^Z $ is defined by 
$ L_0^Z = \mathrm{id} $ and 
\begin{align}\label{LNTF}
L_{n}^Z
=
- \sum_{k=1}^{n} \binom{n}{k} \hat{H}_{n-k}\left( \int_{0}^{1} \!\! \dZ(s) dW(s), \,\,
\Vert Z \Vert_{\mathscr{H}}^{2} \right) D_{-Z}^{\otimes k}, \quad n \in \mathbf{N} 
\end{align}
where the polynomials $\hat{H}_{n}(x,y)$, $n=1,2,\cdots ,$ are defined 
by means of the formula
\begin{align*}
e^{\lambda x-\frac{\lambda ^{2}}{2}y^{2}} = \sum_{n=0}^{\infty} \frac{\lambda^{n}}{n!} \hat{H}_{n}(x,y).
\end{align*}
With this notation, the Hermite polynomials we have used so far are can be written as
\begin{align*}
H_{n}[x] = \hat{H}_{n}(x,1).
\end{align*}

\begin{theorem}\label{new_operator}
For any $ F \in \mathbf{P} $, we have
\begin{align}\label{new_operator_exponential}
&
E \Big[\sum_{n=0}^\infty \frac{1}{n!}
L_{n}^Z F 
\Big] 
= E \Big[ \mathcal{E} ( \int_{0}^{\cdot} \!\! \dZ (s) dW(s) )_{1} \cdot F 
\Big].
\end{align}
\end{theorem}
\begin{proof}
It suffices to show 
\begin{align}\label{new_operator_exponential0}
&
E \Big[
L_{n}^Z F 
\Big] 
= E \Big[ \hat{H}_{n}
\Big( \int_{0}^{1} \!\! \dZ(s) dW_{s}, 
\Vert Z \Vert_{\mathscr{H}}^{2} 
\Big) \cdot F 
\Big]
\end{align}
for each $ n \in \mathbf{N} $ and $ F \in \mathbf{P} $. 
If we can prove that 
\begin{align}\label{new_operator_exponential2}
E
\Big[
L_{n}^Z \Big( \mathcal{E}( \int \! \df dW )_{1} \Big)
\Big]
= E
\Big[ \hat{H}_{n}
\Big( \int_{0}^{1} \!\! \dZ(s) dW_{s}, 
\Vert Z \Vert_{\mathscr{H}}^{2}  \Big)
\cdot \mathcal{E} ( \int \! \df dW )_{1}
\Big]
\end{align}
for arbitrary $ f \in \mathscr{H} $, then (\ref{new_operator_exponential0})
is deduced. In fact, for a finite orthonormal system $ \{e_1, \cdots, e_m \} $, 
take $ f := \lambda_1 e_1 + \cdots \lambda_m e_m $ for $ \lambda_1, \cdots, \lambda_m \in \mathbf{R} $. 
Then, 
\begin{equation*}
\begin{split}
\mathcal{E}( \int \! \df dW )_{1}
&= \prod_{i=1}^m  \mathcal{E}( \lambda_{i} \! \int \! \de_{i} dW )_{1} \\
&= \sum_{N=0}^\infty 
\frac{1}{N!} \sum_{n_1 + \cdots +n_m =N }
\frac{N!}{n_{1}! \cdots n_{m}!}
\prod_{i=1}^m \lambda_{i}^{n_i} H_{n_i}
( \int_{0}^{1} \! \de_{i}(s) dW(s) ),
\end{split}
\end{equation*}
and we notice that $ \sum_{N=0}^\infty a_N  $ where 
\begin{equation*}
a_N
= E \left[ \sum_{n_1 + \cdots +n_m =N }
\frac{N!}{n_{1}! \cdots n_{m}!}
\prod_{i=1}^m \lambda_{i}^{n_i} H_{n_i}
( \int_{0}^{1} \!\! \de_i (s) dW(s) ) \right]
=\left\{\begin{array}{cc}
1 & \text{if } N=0 \\
0 & \text{otherwise}
\end{array}\right.
\end{equation*}
is absolutely convergent. This means that 
(\ref{new_operator_exponential0}) is valid for arbitrary monomials 
and hence for all polynomials. 

So, let us prove (\ref{new_operator_exponential2}). First we note that  
\begin{align*}
&
E \Big[ L_{n}^Z \Big( \mathcal{E}( \int \! \df dW )_{1} \Big) \Big] \\
&=
E \Big[ \sum_{k=1}^{n} (-1)^{k+1} \binom{n}{k} \hat{H}_{n-k} \Big( \int_{0}^{1} \!\!\! \dZ(s) dW_{s},
\Vert Z \Vert_{\mathscr{H}}^{2} \Big) D_{Z}^{\otimes k} \mathcal{E}( \int \! \df dW )_{1} \Big] ,
\end{align*}
where $\hat{H}_{n}(s)$ denotes 
$\hat{H}_{n}(\int_{0}^{s}\dZ(u)dW_{u},\int_{0}^{s}\dZ(u)^{2}du)$ and $\hat{H}_{n}:=\hat{H}_{n}(1)$. 
Since $D_{i}\mathcal{E}( \int \df dW )_{1}= \langle f, e_i\rangle \mathcal{E}( \int \df dW )_{1}$, we have
\begin{align*}
&
E \Big[ L_{n}^{Z} \Big( \mathcal{E}( \int \! \df dW )_{1} \Big) \Big] \\
&= \nonumber
E \Big[ \mathcal{E}( \int \! \df dW )_{1}
\Big\{
\sum_{k=1}^{n}(-1)^{k+1} \binom{n}{k}
\hat{H}_{n-k} \!\!\!\!
\sum_{i_{1},\cdots ,i_{k}} \!\!\!
Z_{i_1} \cdots Z_{i_k} \langle f, e_{i_1} \rangle \cdots \langle f, e_{i_k} \rangle
\Big\} \Big] \\
&= \nonumber
 E \Big[ \mathcal{E}( \int \! \df dW )_{1}
\Big\{
\sum_{k=1}^{n} (-1)^{k+1} \binom{n}{k}
\hat{H}_{n-k} \langle Z, f \rangle ^{k}
\Big\}
\Big]. 
\end{align*}
We will use the following formulas to obtain (\ref{new_operator_exponential2}) which will complete the proof;
\begin{align*}
\hat{H}_{n}(t) = n \int_{0}^{t} \! \hat{H}_{n-1} (s) \dZ(s) dW(s),
\end{align*}
\begin{align*}
\mathcal{E} \Big( \int \!\! \df dW \Big) _{t}
= 1+ \int_{0}^{t} \!\! \mathcal{E} ( \int \!\! \df dW )_{s}\df(s)dW(s),
\end{align*}
and
\begin{align}\label{qdr1}
d\langle \hat{H}_{n}, \mathcal{E}( \int \!\! \df dW ) \rangle _{s}
=
n\hat{H}_{n-1}(s) \mathcal{E}( \int \!\! \df dW )_{s} \df(s) \dZ(s) ds.
\end{align}

As a first step we have
\begin{align*}
& E
\Big[
\hat{H}_{n} ( \int_{0}^{1} \!\! \dZ(s) dW_{s}, 
\int_{0}^{1} \!\! \dZ(s)^{2} ds )
\cdot \mathcal{E} ( \int \! \df dW )_{1}
\Big] \\
& = E \Big[ n \int_{0}^{1} \!\! \hat{H}_{n-1}(s)\dZ (s)dW(s) \Big] \\
&\hspace{10mm}+
E \Big[
n \int_{0}^{1} \!\! \hat{H}_{n-1}(s) \dZ(s)dW(s)
\int_{0}^{1} \!\! \mathcal{E}
( \int \!\! \df dW )_{s} \df(s) dW(s) \Big] \\
& =
E \Big[
n \int_{0}^{1} \!\! \hat{H}_{n-1}(s) \mathcal{E} ( \int \!\! \df dW )_{s} \df(s) \dZ(s) ds
\Big] =: I.
\end{align*}
By Ito's formula, we have
\begin{align*}
&
\hat{H}_{n-1}(1) \mathcal{E} ( \int \! \df dW )_{1} 
\int_{0}^{1} \!\! \df(s) \dZ(s) ds \\
&\hspace{0mm}=
\int_{0}^{1} \!\! \hat{H}_{n-1}(s) \mathcal{E}\Big( \int \!\! \df dW \Big)_{s} 
\df(s) \dZ(s) ds
+
\int_{0}^{1} \!\!\int_{0}^{s} \!\! \df(u) \dZ(u) du \ 
d\langle \hat{H}_{n-1}, \mathcal{E}( \int \!\! \df dW ) \rangle _{s} \\
&\hspace{10mm}+ \text{a martingale}.
\end{align*}
Then by using (\ref{qdr1}), we have
\begin{align*}
I
&=
E \Big[ n\hat{H}_{n-1} \mathcal{E}( \int \!\! \df dW )_{1}
\int_{0}^{1} \!\! \df(s) \dZ(s) ds \Big] \\
&\hspace{10mm}-
E \Big[ n(n-1) \int_{0}^{1} \!\! \df(s) \dZ(s)\int_{0}^{s} \!\! \df(u) \dZ(u)du
\ \hat{H}_{n-2}(s) \ \mathcal{E} ( \int \!\! \df dW )_{s}ds \Big] \\
&=:
E \Big[ n\hat{H}_{n-1} \mathcal{E}( \int \!\! \df dW )_{1} \langle f, Z \rangle \Big]
- I\!\!I .
\end{align*}
Again we apply Ito's formula to get
\begin{align*}
&
\hat{H}_{n-2}(1) \mathcal{E} ( \int \!\! \df dW )_{1} 
\langle f,Z \rangle^{2} \\
&=
2\int_{0}^{1} \!\! \hat{H}_{n-2}(s) \mathcal{E} ( \int \!\! \df dW )_{s} 
\int_{0}^{s} \!\! \df(u) \dZ(u)du
\  f(s)Z(s) ds \\
&\hspace{5mm}+
\int_{0}^{1} \!\! \Big\{ \int_{0}^{s} \!\! \df(u) \dZ(u)du \Big\}^{2}
d\langle \hat{H}_{n-2}, \mathcal{E}( \int \!\! \df dW ) \rangle _{s}
+\text{a martingale}
\end{align*}
and by using (\ref{qdr1}) again, we obtain
\begin{align*}
I\!\!I
&=
E\Big[
\frac{n(n-1)}{2} \hat{H}_{n-2} \mathcal{E} ( \int \!\! \df dW )_{1} 
\langle f,Z \rangle^{2}
\Big] \\
&\hspace{5mm}-
E\Big[
\frac{n(n-1)(n-2)}{2}\int_{0}^{1} \!\! \hat{H}_{n-3}(s)
\mathcal{E}( \int \!\! \df dW )_{s} \df(s) \dZ(s)
\Big\{ \int_{0}^{s} \!\! \df(u) \dZ(u) du \Big\}^{2} ds.
\Big]
\end{align*}
Hence we have
\begin{align*}
&E\Big[
\hat{H}_{n}(\int_{0}^{1} \!\! \dZ (s)dW_{s}, \int_{0}^{1} \!\! \dZ(s)^{2}ds )
\cdot \mathcal{E}( \int \!\! \df dW )_{1}
\Big] = I \\
&=
E \Big[
n\hat{H}_{n-1} \mathcal{E} ( \int \!\! \df dW )_{1}
\langle f, Z \rangle
\Big] \\
&\hspace{5mm}-
E \Big[
\frac{n(n-1)}{2} \hat{H}_{n-2} \mathcal{E} ( \int \!\! \df dW )_{1}  \langle f, Z \rangle ^{2}
\Big] \\
&\hspace{3mm}+
E \Big[
\frac{n(n-1)(n-2)}{2}
\int_{0}^{1} \!\! \df(s) \dZ(s) 
\Big\{ \int_{0}^{s} \!\! \df(u) \dZ(u) du \Big\}^{2}
\ \hat{H}_{n-3} (s) \mathcal{E} ( \int \!\! \df dW )_{s}ds \Big]. \\
\end{align*}
By repeating this procedure until $ \hat{H}_* (s) $ in the integrand vanishes, we obtain
\begin{align*}
&
E \Big[
\hat{H}_{n} \Big( \int_{0}^{1} \!\! Z(s) dW_{s}, \int_{0}^{1} \!\! Z(s)^{2} ds \Big)
\cdot \mathcal{E} ( \int \! \df dW )_{1} \Big] \\
&\hspace{10mm}=
E \Big[ \mathcal{E}( \int \! f dW )_{1}
\Big\{ \sum_{k=1}^{n}(-1)^{k+1} \binom{n}{k}
\hat{H}_{n-k} \langle Z, f \rangle ^{k} \Big\} \Big].
\end{align*}
\end{proof}

\subsection{Passage to the Cameron-Martin-Maruyama-Girsanov formula.}

From Proposition \ref{heikou_vector} and Theorem \ref{new_operator}, we will give a new proof of Maruyama-Girsanov formula in the case of $f \in \mathbf{P}$.

\begin{lemma}\label{cmmgtork}
As an operator acting on $ \mathbf{P} $, 
\begin{equation*}
 \sum_{n=1}^\infty \frac{1}{n!} L_{n}^{Z}
= \exp \Big\{\int _{0}^{1} 
\!\! \dot{Z}(t)dW(t) - \frac{1}{2}\int _{0}^{1} \!\! \dot{Z}(t)^{2}dt
\Big\} (1-\widetilde{e}^{D_{Z}}) .
\end{equation*}
\end{lemma}
\begin{proof}
\begin{align*}
\sum_{n=0}^{\infty} \frac{1}{n!} L_{n}^{Z}
&=
1 - \sum_{n=1}^{\infty} \frac{1}{n!} \sum_{k=1}^{n} \binom{n}{k}
\hat{H}_{n-k} ( \int_{0}^{1} \!\! \dZ (s) dW(s), \int_{0}^{1} \!\! \dZ (s)^{2} ds )
D_{-Z}^{\otimes k} \\
&\hspace{-5mm}=
1 - \sum_{k=1}^{\infty} \Big( \sum_{n=k}^{\infty} \frac{1}{k!(n-k)!}
\hat{H}_{n-k} ( \int_{0}^{1} \!\! \dZ (s) dW(s), \int_{0}^{1} \!\! \dZ (s)^{2} ds ) \Big)
D_{-Z}^{\otimes k} \\
&\hspace{-5mm}=
1 - \sum_{k=1}^{\infty} \frac{1}{k!}
\Big( \sum_{m=0}^{\infty} \frac{1}{m!}
\hat{H}_{m} ( \int_{0}^{1} \!\! \dZ (s) dW(s), \int_{0}^{1} \!\! \dZ (s)^{2} ds ) \Big)
D_{-Z}^{\otimes k} \\
&\hspace{-5mm}=
1 - \mathcal{E} ( \int \!\! \dZ dW )_{1}
\sum_{k=1}^{\infty} \frac{1}{k!} D_{-Z}^{\otimes k} \\
&\hspace{-5mm}=
1- \mathcal{E} ( \int \!\! \dZ dW )_{1}
\sum_{k=0}^{\infty} \frac{1}{k!} D_{-Z}^{\otimes k}
+ \mathcal{E} ( \int \!\! \dZ dW )_{1}. 
\end{align*}
\end{proof}

\begin{corollary}[Cameron-Martin-Maruyama-Girsanov formula]\label{c_m_m_g}
For $f \in \mathbf{P}$, the following formula holds
\begin{equation}\label{m_g_theorem_vector_field}
E \Big[
\mathcal{E} ( \int \!\! \dZ dW )_{1} f \Big( W-\int_{0}^{\cdot} \!\! \dZ(s) ds \Big)
\Big]
= E \Big[ f(W) \Big].
\end{equation}
\end{corollary}
\begin{proof}
By Lemma \ref{cmmgtork}, we have
\begin{align}\label{vector_field_CMMG2}
&
E \Big[ \sum_{n=0}^{\infty} \frac{1}{n!} L_{n} \Big( f(W) \Big) \Big] \\
&=
E \Big[ f(W) - \mathcal{E} ( \int \!\! \dZ dW )_{1}
\sum_{k=0}^{\infty} \frac{1}{k!} D_{-Z}^{\otimes k} f(W)
+ \mathcal{E} ( \int \!\! \dZ dW )_{1} f(W) \Big] \nonumber \\
&=
E \Big[ f(W) - \mathcal{E} ( \int \!\! \dZ dW )_{1} \widetilde{e}^{D_{-Z}}f(W)
+ \mathcal{E} ( \int \!\! \dZ dW )_{1} f(W) \Big] \nonumber \\
&=
E \Big[ f(W) - \mathcal{E} ( \int \!\! \dZ dW )_{1}
f \Big( W-\int_{0}^{\cdot} \!\! \dZ(s) ds \Big)
+ \mathcal{E} ( \int \!\! \dZ dW )_{1} f(W) \Big]. \nonumber
\end{align}
Then by Theorem \ref{new_operator}, we obtain (\ref{m_g_theorem_vector_field}).
\end{proof}

\section{Another Algebraic Proof for CMMG Formula.}\label{girsanov}

As we have mentioned in the introduction, we give an alternative proof 
which is ``purely" algebraic in the sense that we do not use 
stochastic calculus essentially, though we restrict ourselves in the case
of piecewise constant (=finite-dimensional) case. 

Let $ \mathcal{F} \equiv \{ \mathcal{F}_t \}_{0 \leq t \leq 1} $ be 
the natural filtration of $\mathscr{W}$. 
Let us consider a simple $ \mathcal{F} $-predictable 
process 
\begin{equation}\label{SPP}
z (w,t) = \sum _{k=1}^{2^{s}} 2^{s/2} z_{k}(w) \, 1_{( \frac{k-1}{2^{s}} ,\frac{k}{2^{s}}] }(t)
\end{equation}
where $ z_k $, $ k=1, \cdots, 2^s $ are 
$ \mathcal{F}_{\frac{k-1}{2^{s}}} $-
measurable random variables. 
Define $ \sigma^s_k \in \mathscr{H} $, $ k=1, \cdots, 2^s $ by
\begin{equation*}
\sigma^s_k (t) := 2^{s/2} \int_0^t 
1_{ ( \frac{k-1}{2^{s}} ,\frac{k}{2^{s}}] }(u) \,du.
\end{equation*}
We will suppress the superscript $ s $ whenever it is clear from the context. 
Clearly, 
\begin{equation}\label{cmb}
D_{\sigma_k} F = 0 
\end{equation}
for any $ \mathcal{F}_{\frac{k-1}{2^{s}}} $-measurable
random variable $F$. 
Put
\begin{equation*}
D_{z_k} := z_k D_{\sigma_k} \ \text{and}\  
D^*_{z_k} := z_k D^*_{\sigma_k},
\end{equation*}
for $ k=1, \cdots, 2^{s} $.  
 
\begin{lemma}\label{adj001}
For any $ n \in \Natural $ and $ f \in \mathbf{P} $, we have
\begin{equation}\label{PVF1}
D_{z_{k}}^{n} f = \underbrace{z_{k} D_{\sigma_k} 
\cdots z_{k} D_{\sigma_k}}_{n\ \text{times}} f 
= z_{k}^{n} D_{\sigma_k}^{n} f
\end{equation}
and 
\begin{equation}\label{PVF2}
(D_{z_{k}}^*)^{n} f = \underbrace{z_{k} D^*_{\sigma_k} 
\cdots z_{k} D^*_{\sigma_k}}_{n\ \text{times}} f 
= z_{k}^{n} (D^*_{\sigma_k})^{n} f
\end{equation}
\end{lemma}
\begin{proof}
These are direct from the following ``commutativity":
\begin{equation*}
D_{\sigma_j} (z_i f) = z_i D_{\sigma_j} f,\,\, \text{and}\,\,
D^*_{\sigma_j} (z_i f) = z_i D^*_{\sigma_j} f, \quad \text{if $ i \leq j $}
\end{equation*}
for differentiable $ f $. These follows 
since $ D_{\sigma_j} (z_i) = 0 $.  
\end{proof}

Define the exponentials as 
\begin{equation*}
e^{D_{z_{k}}}
:=\sum_{n=0}^{\infty} \frac{1}{n!} 
D_{z_{k}}^{n}
\quad k=1,2,\cdots ,N
\end{equation*}
and 
\begin{equation*}
e^{D_{z_{k}}^{*}}:=\sum _{n=0}^{\infty} \frac{1}{n!} 
(D_{z_{k}}^*)^n,\quad k=1,2,\cdots N
\end{equation*}
formally. 
By Lemma \ref{adj001} we have 
\begin{equation*}
e^{D_{z_{k}}}=\sum_{n=0}^{\infty} \frac{z^n_k}{n!} 
D_{\sigma_{k}}^{n}
\end{equation*}
and thus we can include $ \mathbf{P} $ in the domain
of $ e^{D_{z_{k}}} $.  

Let us introduce a subspace $ \mathbf{P}_H $ of $ \mathbf{P} $, 
which consists of polynomials with respect to 
$ \{ [e_i] (w)\} $, where $\{ e_i \} $ is the Haar system.
Note that $ \mathbf{P}_H $ is 
also characterized as 
all the polynomials with respect to 
$ \{ [\dot{\sigma}^s_k] (w): k=1,\cdots,2^s, s \in \Natural \} $. 

The following is a main result in our program.
\begin{theorem}\label{adj007}
(i) For any $ F \in \mathbf{P}_H $, we have
\begin{equation}\label{Taylor01}
e^{D_{z_{2^s}}} \cdots e^{D_{z_{1}}} F(w) 
=F(w + \int_0^\cdot z(w,u)\,du).
\end{equation}
(ii) For any $ \mathcal{F}_{(k-1)/2^s} $-measurable 
random variable $F$,
\begin{equation}\label{commute01}
e^{D_{z_k}^*} F = F e^{D_{z_k}^*} (1).
\end{equation}  
In particular, the function $F$ is in the domain of $ e^{D_{z_k}^*} $. 
Furthermore, we have
\begin{equation}\label{expo01}
e^{D_{z_{2^s}}^{*}}\cdots e^{D_{z_{1}}^{*}}(1) 
= \exp \Big\{ 
\int_{0}^{1}
 z(w,s)dw(s) 
- \frac{1}{2} \int_{0}^{1}
z (w,s)^{2}ds \Big\},
\end{equation}
(iii) Fix $ k \in \Natural $. 
Let $F\in \mathbf{P}$ and let $G$ be an arbitrary $ \mathcal{F}_{(k-1)/2^{s}} $-measurable integrable function.
Then 
\begin{equation}\label{adj006}
E[ e^{D_{z_k}} (F) G ] = E[ F e^{D^*_{z_k}} (G) ]. 
\end{equation}
\end{theorem}
\begin{proof}
(i) First, notice that 
$ F \in \mathbf{P}_H $ 
is always expressed as a linear combination of
$
\prod_{k=1}^{2^s} F_{k},
$
where each $ F_{k} $ is a polynomial in 
\begin{align}
\left\{ [\sigma^t_l] (w):  \Big(\frac{l-1}{2^t}, \frac{l}{2^t} \Big] 
\subset \Big(\frac{k-1}{2^s}, \frac{k}{2^s} \Big] \right\},
\label{haar}
\end{align}
so that we can assume that $F$ is of the form
\begin{equation*}
F = \sum_{i=1}^{N} \prod_{k=1}^{2^s} F_{k,i},
\end{equation*}
where each $ F_{k,i} $ is a polynomial in (\ref{haar}).
By Proposition \ref{Prop1} and the definition of $ D_{\sigma_k} $,
we have 
\begin{equation*}
e^{D_{z_{k}}} F_{l,i} (w)
= \begin{cases}
F_{k,i} (w + z_k \sigma_k) & (l=k) \\
F_{l,i} (w) & (l \ne k).
\end{cases}
\end{equation*}
Then by Lemma \ref{Lem1}, 
\begin{equation*}
e^{D_{z_{k}}} \prod_{l=1}^{2^s} F_{l,i} (w)
= F_{k,i} (w + z_k \sigma_k) \prod_{l\ne k} F_{l,i} (w).
\end{equation*}
Since $ z_k $ is $ \mathcal{F}_{t_k} $-measurable, 
we also have, if $ j > k $, 
\begin{equation*}
\begin{split}
& e^{D_{z_{j}}} e^{D_{z_{k}}} \prod_{l=1}^{2^s} F_{l,i} (w) \\
& = e^{D_{z_{j}}} F_{k,i} (w + z_k \sigma_k) 
e^{D_{z_{j}}} \prod_{l\ne k} F_{l,i} (w) \\
&= F_{k,i} (w + z_k \sigma_k) F_{j,i} (w + z_j \sigma_j )
\prod_{l\ne j,k} F_{l,i} (w).
\end{split}
\end{equation*}
Then, inductively we have 
\begin{equation*}
e^{D_{z_{2^s}}} \cdots e^{D_{z_{1}}} \prod_{l=1}^{2^s} F_{l,i} (w)
= \prod_{l=1}^{2^s} F_{l,i} (w + z_l \sigma_l ),
\end{equation*}
and by linearity we obtain (\ref{Taylor01})
since 
\begin{equation*}
\sum_{l=1}^{2^s} z_{l}(w) \sigma_{l} (t) = \int_0^t z (w,u)\,du. 
\end{equation*}

(ii) Noting that $ D_{\sigma_k} F = 0 $ for $ \mathcal{F}_{(k-1)/2^s} $ -
measurable random variable $ F $, we have
\begin{equation*}
\begin{split}
D^*_{z_k} F &= z_k \{-D_{\sigma_k} +  2^{s/2}(w_{k/2^s}-w_{(k-1)/2^s})\}F \\
&= F  z_k 2^{s/2}(w_{k/2^s}-w_{(k-1)/2^s}) = F D^*_{z_k} (1) 
\end{split}
\end{equation*}
since $ z_k $ is also $ \mathcal{F}_{(k-1)/2^{s}} $-measurable.
Inductively, we then have 
\begin{equation*}
(D^*_{z_k})^n F = F (D^*_{z_k})^n (1),
\end{equation*}
and hence we have (\ref{commute01}), 
which in turn implies (\ref{expo01}).
In fact, we have by induction
\begin{equation*}
e^{D_{z_{2^s}}^{*}}\cdots e^{D_{z_{1}}^{*}}(1) 
= \prod_{k=1}^{2^s} \{ e^{D_{z_{k}}^{*}}(1) \} 
\end{equation*} 
since $ e^{D_{z_{k-1}}^*} \cdots e^{D_{z_{1}}^*} (1) $
is $ \mathcal{F}_{(k-1)/2^s} $-measurable for any $ k $,
and for each $i=1,2,\cdots ,2^s $, we have
\begin{equation*}
\begin{split}
e^{D_{z_{i}}^{*}}(1)&=\sum _{n=0}^{\infty} 
\frac{z_{i}^{n}}{n!}
(D_{\sigma_i}^*)^n (1) = \sum _{n=0}^{\infty} 
\frac{z_{i}^{n}}{n!}H_{n}[\int_0^1 \sigma_k (t) \,dw_t ] \\
& =\exp \Big\{ z_{i}(w) 2^{s/2} (w_{k/2^s}-w_{(k-1)/2^s}) 
-\frac{1}{2}z_{i}(w)^{2} \Big\}.
\end{split}
\end{equation*}

(iii) Since $ F $ is a polynomial, 
\begin{equation*}
e^{D_{z_k}} F
= \sum_{n=0}^M \frac{z_k^n}{n!} D_{\sigma_k}^{n} F 
\end{equation*}
for some $ M \in \Natural \cup \{ 0\} $. 
Therefore, the left-hand-side of (\ref{adj006}) is rewritten as 
\begin{equation*}
\sum_{n=0}^M \frac{1}{n!} 
E[ z_k^n D_{\sigma_k}^n F \cdot G ].
\end{equation*}
Since $ z_k $ and $ G $ are $ \mathcal{F}_{(k-1)/2^s} $-measurable, 
we have, for $ n \leq M $
\begin{equation*}
\begin{split}
& 
E[ z_k^n D_{\sigma_k}^n F \cdot G] 
= E [ F \cdot (D^*_{\sigma_k})^n  
z^n_k G] \\
&= E [ F \cdot z^n_k (D^{*}_{\sigma_k})^n G ]
= E[ F \cdot (D_{z_k}^{*})^n G].
\end{split}
\end{equation*}
The relation is valid for $ n > M $ since
\begin{equation*}
(D_{\sigma_k}^*)^{n} G 
= G (D^{*}_{\sigma_k})^n (1) =  G H_n (\int_0^1 \sigma_k (t) \,dw_t), 
\end{equation*}
and the degree of $ F $ as a polynomial of 
$ \int_0^1 \sigma_k (t) \,dw_t $ is less than $ M $,
we have
\begin{equation*}
E[ z_k^n D_{\sigma_{k}}^{n} F \cdot G]
= E[ F \cdot  D_{z_{k}}^{*n}G ] = 0. 
\end{equation*}
Thus we have
\begin{equation*}
E[ \sum_{n=0}^\infty \frac{1}{n!} D_{z_k}^n F \cdot G]
= E[ \sum_{n=0}^\infty \frac{1}{n!}  F
\cdot D^{*n}_{z_k} G ],
\end{equation*}
which is the desired relation. 
\end{proof}

\begin{remark}
(i) We do not assume smoothness for $ F $ in (\ref{commute01}).
(ii) In (\ref{Taylor01}) and (\ref{expo01}), 
the order of application of the operators is important. 
If it is changed anywhere, neither holds anymore. 
\end{remark}

By using the above algebraic results, we can prove the following
\begin{corollary}[Cameron-Martin-Maruyama-Girsanov formula]\label{MGFP}
For a simple 
predictable $ z $ in (\ref{SPP}) 
and $ F \in \mathbf{P}_H $, it holds
\begin{equation}\label{FMG}
\begin{split}
& E[ F( w -\int_0^\cdot z(w,u)\,du)
\exp \Big\{ \int_0^1 z(w,t)\,dw_t 
-\frac{1}{2} \int_0^1 | z(w,t) |^{2}\,dt \Big\}\\
&\hspace{1cm} = E[F]. 
\end{split}
\end{equation}
\end{corollary}
\begin{proof}
As a formal series, we have
\begin{equation*}
e^{D_{z_k}} e^{-D_{z_k}} =1,
\end{equation*}
for $ k=1,\cdots 2^s $.
Then, for $ F \in \mathbf{P}_H $, we have
\begin{equation*}
F  = e^{D_{z_1}} e^{-D_{z_1}} F
\end{equation*}
and since $ e^{-D_{z_1}} F $ is a polynomial, 
by Theorem \ref{adj007} (iii), we have
\begin{equation}\label{adj009}
\begin{split}
& E[ F ] = E[ e^{D_{z_1}} 
e^{-D_{z_1}} F] \\
&= E[ e^{-D_{z_1}} F \cdot 
e^{D^*_{z_1}} (1) ]. 
\end{split}
\end{equation} 
Inductively, since 
\begin{equation*}
e^{-\partial_{z_{k}}} \cdots e^{-\partial_{z_1}} f (\xi)
\end{equation*}
still is a polynomial in
\begin{equation*}
\left\{ [\sigma^t_l] (w):  \Big(\frac{l-1}{2^t}, \frac{l}{2^t} \Big] 
\subset \Big(\frac{k-1}{2^s}, \frac{k}{2^s} \Big] \right\}, 
\end{equation*}
and 
\begin{equation*}
e^{D^*_{z_{k-1}}} \cdots e^{D^*_{z_1}} (1)
\end{equation*}
is $ \mathcal{F}_{(k-1)/2^s} $-measurable,
we have
\begin{equation}\label{adj010}
\begin{split}
& E[ F ] \\
&= E[ e^{D_{z_k}} e^{-D_{z_k}} 
e^{-D_{z_{k-1}}} \cdots e^{-D_{z_1}} F
\cdot e^{D^*_{z_{k-1}}} \cdots e^{D^*_{z_1}} (1) ] \\
&= E[
e^{-D_{z_{k}}} \cdots e^{-D_{z_1}} F \cdot 
e^{D^*_{z_{k}}} \cdots 
e^{D^*_{z_1}} (1)]. 
\end{split}
\end{equation} 
Combining this with (\ref{Taylor01}) and (\ref{expo01})
in Theorem \ref{commute01}, we have the formula (\ref{FMG}).
\end{proof}


\appendix \section{Continuity of the translation}\label{AP}

The following lemma extends the translation 
on the dense subset of polynomials to an operator on 
$ L_q $ to $ L_p $, and hence ensure the 
MG formula (\ref{FMG}) for any bounded measurable $ F $.

\begin{lemma}
Let $ z $ be a predictable process as (\ref{SPP}).
Suppose that 
\begin{equation}\label{cnovi}
E \left[ \exp \left\{ c \int_{0}^{1} \!\! z(t)^{2} \,dt \right\} \right] < \infty
\end{equation}
for some $ c > 0 $.  
Then, for $ p \in [1, \infty) $, there exists $ q \in (p, \infty) $ 
and a positive constant $ C_p $ such that 
\begin{equation*}
\Vert  e^{-D_{z_{2^s}}} \cdots e^{-D_{z_{1}}} F \Vert_p  
\leq C_p \Vert F \Vert_q 
\end{equation*}
for any $ F \in \mathbf{P}_H $. 
\end{lemma}
\begin{proof}
We will denote $ Z := \int_0^\cdot z(t) \,dt $ and 
\begin{equation*}
\mathcal{E} (z)
:= \exp \left\{\int_0^1 \!\! z(t) \,dw(t) - \frac{1}{2}  \int_0^1 \!\! z(t)^{2} \,dt 
\right\}. 
\end{equation*}
Let $ n \geq 1 $ be an integer and $ p < 2n $. 
By H\"older's inequality, 
\begin{equation*}
\begin{split}
& E \left[ \left| F \left(w - Z (w) \right) \right|^p \right] 
= E \left[ \left| F \left(w - Z(w)\right) \right|^p 
\{\mathcal{E} (z)\}^{\frac{p}{2n}} 
\{\mathcal{E} (z)\}^{-\frac{p}{2n}} \right] \\
& \leq E \left[ \left| F (w - Z(w) ) \right|^{p \cdot \frac{2n}{p}} 
\{\mathcal{E} (z)\}^{\frac{p}{2n} \cdot \frac{2n}{p}} 
\right]^{\frac{p}{2n}} 
\cdot E \left[ \{\mathcal{E} (z)\}^{-\frac{p}{2n} 
\cdot \frac{2n}{2n-p}}
\right]^{\frac{2n-p}{2n}} \\
&= E \left[ | F (w - Z(w) ) |^{2n} 
\mathcal{E} (z) \right]^{\frac{p}{2n}}
\cdot  E \left[ \{\mathcal{E} (z)\}^{-\frac{p}{2n-p}} 
\right]^{\frac{2n-p}{2n}}.
\end{split}
\end{equation*}
Since $ F $ is a polynomial, so is $ |F|^{2n} $. 
Therefore, we can apply the MG formula for polynomials 
(\ref{FMG}) in Corollary \ref{MGFP},  to obtain
\begin{equation*}
E \left[ | F (w - Z(w) ) |^{2n} 
\mathcal{E} (z) \right]^{\frac{p}{2n}}
= E \left[ |F|^{2n} \right]^{\frac{p}{2n}}
= \Vert F \Vert_{2n}^p. 
\end{equation*}
Now it suffices to show that 
\begin{equation}\label{cest}
E \left[ \{\mathcal{E} (z)\}^{-\frac{p}{2n-p}} 
\right] < \infty. 
\end{equation}
Let us denote $ L_t:= \int_0^t z(u)\, dw(u) $. Then
$ \langle L \rangle_t = \int_0^t z(u)^{2} \,du $. 
Now, since we have 
\begin{equation*}
\begin{split}
& \{\mathcal{E} (z)\}^{-\frac{p}{2n-p}} 
= \exp \left\{
-\frac{p}{2n-p} L - \frac{p^2}{(2n-p)^2} \langle L \rangle 
\right\} \\
& \hspace{2cm} \cdot
\exp \left\{ \left(\frac{p}{2(2n-p)} + \frac{p^2}{(2n-p)^2}
\right)  \langle L \rangle \right\}, 
\end{split}
\end{equation*}
by Schwartz inequality we have
\begin{equation*}
\begin{split}
& E \left[ \{\mathcal{E} (z)\}^{-\frac{p}{2n-p}} 
\right] \\ & \leq 
E \left[\exp \left\{
-\frac{2p}{2n-p} L - \frac{2p^2}{(2n-p)^2} \langle L \rangle 
\right\}
\right]^{1/2} \\
& \hspace{1cm} 
\cdot E \left[ \exp 
\left\{ \left(\frac{p}{(2n-p)} + \frac{2p^2}{(2n-p)^2}
\right)  \langle L \rangle \right\}
\right]^{1/2}. 
\end{split}
\end{equation*}
Clearly, $ \frac{p}{(2n-p)} + \frac{2p^2}{(2n-p)^2} \to 0 $ as 
$ n \to \infty $, and hence we can take large enough $ n $ to have 
the estimate (\ref{cest}) by using the assumption (\ref{cnovi}). 
\end{proof}
\begin{remark}
By a similar but easier procedure we can also prove 
a continuity lemma for $ e^{D_\theta} $ with
$ \theta \in \mathscr{H} $, to extend (\ref{eq:exp2}) 
in Corollary \ref{Cor1}
to obtain a full version of CM formula.  
\end{remark}

\end{document}